\documentclass[a4paper,11pt]{article}
\usepackage{amsmath,amsthm,amssymb,enumerate}
\usepackage{graphicx}

\numberwithin{equation}{section}
\numberwithin{figure}{section}
\newtheorem{thm}{Theorem}[section]
\newtheorem{lem}[thm]{Lemma}

\newtheorem{definition}[thm]{{Definition}}
\newtheorem{remark}[thm]{{Remark}}
\newtheorem{example}[thm]{\rm{Example}}

\newenvironment{df}{\begin{definition}\rm}{\end{definition}}

\newcommand{\Proof}{\noindent {\it Proof}. \quad}

\title{Triple linking numbers and triple point numbers 
of certain $T^2$-links}
\author{Inasa Nakamura
\\
\\
{\small Research Institute for Mathematical Sciences, Kyoto University,} \\
{\small Kyoto 606-8502, Japan} \\
{\small inasa@kurims.kyoto-u.ac.jp}
}
\date{}

\begin{document}

\maketitle

\begin{abstract}
The triple linking number of an oriented surface link was defined as an analogical notion of the linking number of a classical link. We consider a certain $m$-component $T^2$-link ($m \geq 3$) determined from two commutative pure $m$-braids $a$ and $b$. We present the triple linking number of such a $T^2$-link, by using the linking numbers of the closures of $a$ and $b$. This gives a lower bound of the triple point number. In some cases, we can determine the triple point numbers, each of which is a multiple of four. 
\end{abstract}

%%%%%%%%%%%%%%%%%%%%%%%%%%
\section{Introduction} \label{intro}
%%%%%%%%%%%%%%%%%%%%%%%%%%%
A {\it surface link} is a smooth embedding of a closed surface into the Euclidean 4-space $\mathbb{R}^4$. A {\it $T^2$-link} is a surface link each of whose components is of genus one. In this paper, we consider a certain $m$-component $T^2$-link which is determined from two commutative pure $m$-braids $a$ and $b$. 
The triple linking number of an oriented surface link is defined in \cite{CJKLS} 
as an analogical notion of the linking number of a classical link. The aim of this paper is to present the triple linking number of such a $T^2$-link, by using the linking numbers of the closures of $a$ and $b$. Further, we study the triple point number. The triple linking numbers give a lower bound of the triple point number. In some cases, we can determine the triple point numbers, each of which is a multiple of four. 

\begin{figure}
\begin{center}
 \includegraphics*{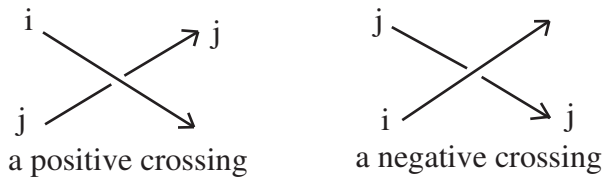}%%%%%%%%%%%%%%Fig.2.5%%%%%%%%%%%%%%%
\end{center}
 \caption{A positive crossing and a negative crossing. }
 \label{0104-2}
  \end{figure}

We review the linking number of an oriented classical link $L$ as follows. 
For positive integers $i$ and $j$ with $i \neq j$, the {\it linking number} of the $i$th and $j$th components of $L$ is the total number of positive crossings minus the total number of negative crossings of a diagram of $L$ such that the over-arc (resp. under-arc) is from the $i$th (resp. $j$th) component; see Fig. \ref{0104-2}. We denote it by $\mathrm{Lk}_{i,j}(L)$. It is known \cite{Rolfsen} that $\mathrm{Lk}_{j,i}(L)=\mathrm{Lk}_{i,j}(L)$. 

\begin{figure}
\begin{center}
  \includegraphics*{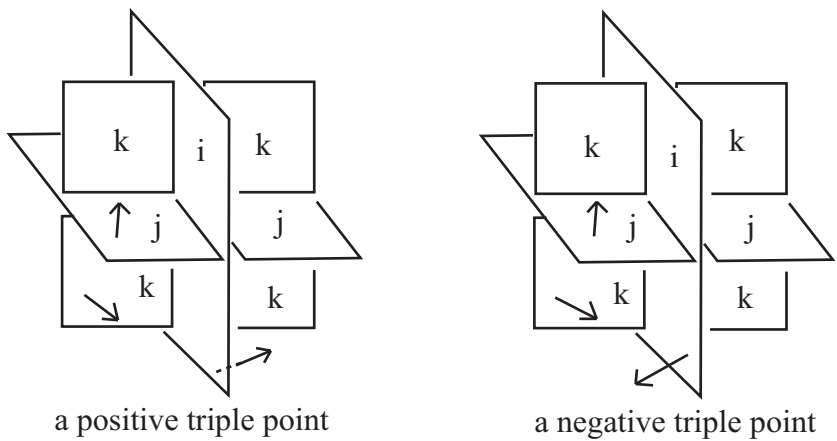}%%%%%%%%%%%%%%Fig.2.5%%%%%%%%%%%%%%%
\end{center}
  \caption{A positive triple point and a negative triple point, where we denote the orientations of sheets by normals.}
  \label{0104-1}
  \end{figure}

The triple linking number of an oriented surface link $S$ is defined as follows (see \cite[Definition 9.1]{CJKLS}, see also \cite{CKS}). 
For positive integers $i$, $j$, and $k$ with $i \neq j$ and $j \neq k$, the {\it triple linking number} of the $i$th, $j$th, and $k$th components of $S$ is the total number of positive triple points minus the total number of negative triple points of a surface diagram of $S$ such that the top, middle, and bottom sheet are from the $i$th, $j$th, and $k$th component of $S$ respectively (\cite{CJKLS}); see Fig. \ref{0104-1}. 
We denote it by $\mathrm{Tlk}_{i,j,k}(S)$. 
It is known \cite{CJKLS} that $\mathrm{Tlk}_{k,j,i}(S)=-\mathrm{Tlk}_{i,j,k}(S)$ if $i$, $j$, $k$ are mutually distinct, and otherwise $\mathrm{Tlk}_{i,j,k}(S)=0$. 
 
We show the following theorem. We consider a torus-covering $T^2$-link, which is a $T^2$-link in the form of an unbranched covering over the standard torus (see Definition \ref{Def2-1}). It is known \cite[Lemma 2.8]{N} that a torus-covering $T^2$-link is determined from two commutative classical $m$-braids, which we call basis $m$-braids (see Section \ref{0913-1}). We denote by $\mathcal{S}_m(a,b)$ the torus-covering $T^2$-link with basis $m$-braids $a$ and $b$. In this paper, we consider the case when the basis braids are pure $m$-braids for $m \geq 3$. For an $m$-braid $c$, let us denote by $\hat{c}$ the closure of $c$. 

\begin{thm} \label{linking}
Let $a$ and $b$ be commutative pure $m$-braids for $m \geq 3$. Then the triple linking number $Tlk_{i,j,k}(\mathcal{S}_m(a,b))$ ($i \neq j$ and $j \neq k$) is given by 
\begin{equation*} 
Tlk_{i,j,k}(\mathcal{S}_m(a, b))=-\mathrm{Lk}_{i,j}(\hat{a}) \mathrm{Lk}_{j,k}(\hat{b}) +\mathrm{Lk}_{i,j}(\hat{b}) \mathrm{Lk}_{j,k}(\hat{a}), 
\end{equation*}
where $\mathrm{Lk}_{i,j}(\hat{a})$ (resp. $\mathrm{Lk}_{i,j}(\hat{b})$) is the linking number of the $i$th and $j$th components of $\hat{a}$ (resp. $\hat{b}$). Here we define the $l$th component of $\mathcal{S}_m(a, b)$ (resp. $\hat{c}$ for $c=a$ or $b$) by the component containing the $l$th string of the basis braids (resp. $c$) for $l=1,2,\ldots,m$. 
\end{thm}

The {\it triple point number} of a surface link $S$, denoted by $t(S)$, is the minimal number of triple points among all possible generic projections of $S$. For the results known about triple point numbers, see \cite[Section 4.4.2]{CKS}. By definition, we can see that $t(S) \geq \sum_{i \neq j, j \neq k}|\mathrm{Tlk}_{i,j,k}(S)|$; thus the triple linking numbers of $S$ give a lower bound of the triple point number of $S$. 
In particular, we show the following theorem. 

\begin{thm} \label{thm-triple-pt}%%%%%%%%%%%%%%%
Let $m \geq 3$. Let $b$ be a pure $m$-braid, let $\Delta$ be a full twist of a bundle of $m$ parallel strings, and let $n$ be a non-negative integer. Put $\mu=\sum_{i<j}|\mathrm{Lk}_{i,j}(\hat{b})|$, and put $\nu=\sum_{i<j<k} (\nu_{i,j,k}+\nu_{j,k,i}+\nu_{k,i,j})$, where $\nu_{i,j,k}=\min_{i,j,k} \{ |\mathrm{Lk}_{i,j}(\hat{b})|, |\mathrm{Lk}_{j,k}(\hat{b})| \}$ if $\mathrm{Lk}_{i,j}(\hat{b}) \mathrm{Lk}_{j,k}(\hat{b})>0$ and otherwise zero. Then  
\[
t(\mathcal{S}_m(b, \Delta^{n})) \geq 4n(\mu(m-2)-\nu). 
\]
\end{thm}%%%%%%%%%%%%%%%%%

In some cases, we can determine the triple point numbers. Let $\sigma_1, \sigma_2, \ldots, \sigma_{m-1}$ be the standard generators of the $m$-braid group. 
\begin{thm} \label{example}
Let $m \geq 3$. Let $b$ be an $m$-braid presented by a braid word which consists of $\sigma_i^{2(-1)^i}$ ($i=1,2,\ldots, m-1$); note that $b$ is a pure braid. Then 
\[
t(\mathcal{S}_m(b, \Delta^{n}))=4n(m-2)(\sum_{i<j}|\mathrm{Lk}_{i,j}(\hat{b})|). 
\]
Further the triple point number is realized by a surface diagram in the form of a covering over the torus. 
\end{thm}
It is known \cite{CKSS01} (see also \cite{CKSS02}) that any oriented surface link is bordant to the split union of oriented \lq\lq necklaces", and \cite{CKSS02} any surface link is unorientedly bordant to the split union of necklaces and connected sums of standard projective planes; see also \cite{Sanderson}. A necklace has the triple point number $4n$ (see Section \ref{0104-3}). For other examples of surface links (not necessarily orientable) which realize large triple point numbers, see \cite{COS, Kamada-Oshiro, Oshiro10, Satoh}. In the papers \cite{COS, Kamada-Oshiro, Oshiro10} (resp. \cite{Satoh}), they use quandle cocycle invariants (resp. normal Euler numbers) to give lower bounds of triple point numbers. Quandle cocycle invariants \cite{CJKLS, CKS01, CKS} can be regarded as an extended notion of triple linking numbers (\cite{CJKLS, CKSS01}), useful to give lower bounds of triple point numbers; see \cite{COS, Hatakenaka, Iwakiri, Kamada-Oshiro, Oshiro10, SatShi04, SatShi05}. 

This paper is organized as follows. In Section \ref{section1}, we give the definition of a torus-covering $T^2$-link. In Section \ref{TripleLinking}, we prove Theorem \ref{linking}. In Section \ref{0104-3}, we prove Theorems \ref{thm-triple-pt} and \ref{example}. 

%%%%%%%%%%%%%%%%%%%%%%%%%%%%%%%%%%%%%%%%
\section{A torus-covering $T^2$-link} \label{section1}%%%%%%%%%%%%%%%%%%%%%
%%%%%%%%%%%%%%%%%%%%%%%%%%%%%%%%%%%%%%% 
%
It is known \cite{Kamada1, Kamada3} that any oriented surface link can be presented by a branched covering over the standard 2-sphere. A torus-covering link was introduced in \cite{N} as a new construction of a surface link, by considering the standard torus instead of the standard 2-sphere. In this section, we define a torus-covering $T^2$-link (Definition \ref{Def2-1}, see \cite{N}). The triple linking number of the $i$th, $j$th, and $k$th components of an oriented surface link is the sum of the signs of all the triple points of type $(i,j,k)$ (see Section \ref{0203-1}, see also \cite{CJKLS}). We can obtain the sign and type of each triple point of a surface diagram of $\mathcal{S}_m(a,b)$, by studying the braid word transformation sequence from $ab$ to $ba$ (Lemma \ref{0125-2}). 

%%%%%%%%%%%%%%%%%%%%%%%%%%%%%%%%%%%%%%%%
\subsection{A torus-covering $T^2$-link} \label{0913-1}%%%%%%%%%%%%%%%%%%%%%%
%%%%%%%%%%%%%%%%%%%%%%%%%%%%%%%%%%%%
Let $T$ be the standard torus in $\mathbb{R}^4$, i.e. the boundary of the standard solid torus in $\mathbb{R}^3 \times \{0\}$. Let us fix a point $x_0$ of $T$, and take a meridian $\mathbf{m}$ and a longitude $\mathbf{l}$ of $T$ with the base point $x_0$. 
A {\it meridian} is an oriented simple closed curve on $T$ which bounds a 2-disk in the solid torus whose boundary is $T$ and which is not null-homologous in $T$. A {\it longitude} is an oriented simple closed curve on $T$ which is null-homologous in the complement of the solid torus in the three space $\mathbb{R}^3 \times \{0\}$ and which is not null-homologous in $T$. Let $N(T)$ be a tubular neighborhood of $T$ in $\mathbb{R}^4$. 

\begin{df} \label{Def2-1} %%%%
A {\it torus-covering $T^2$-link} is a surface link $S$ in $\mathbb{R}^4$ such that $S$ is embedded in $N(T)$ and $p |_{S} \,:\, S \rightarrow T$ is an unbranched covering map, where $p \,:\, N(T) \rightarrow T$ is the natural projection. 
\end{df}%%%%%%%%%%%%%%%%%%%
For the definition and properties of a torus-covering link whose component might be of genus more than one, see \cite{N}. 
\\

Let us consider a torus-covering $T^2$-link $S$. The intersections $S \cap p^{-1}(\mathbf{m})$ and $S \cap p^{-1}(\mathbf{l})$ are closures of classical braids. Cutting open the solid tori $p^{-1}(\mathbf{m})$ and $p^{-1}(\mathbf{l})$ at the 2-disk $p^{-1}(x_0)$, we obtain a pair of classical braids. We call them {\it basis braids} (\cite{N}). The basis braids of a torus-covering $T^2$-link are commutative, and for any commutative $m$-braids $a$ and $b$, there exists a unique torus-covering $T^2$-link with basis braids $a$ and $b$ (\cite[Lemma 2.8]{N}). For commutative $m$-braids $a$ and $b$, we denote by $\mathcal{S}_m(a,b)$ the torus-covering $T^2$-link with basis $m$-braids $a$ and $b$ (\cite{N}). 

In this paper, we consider the case when the basis braids are pure braids, and we define the $l$th component of $\mathcal{S}_m(a, b)$ by the component containing the $l$th string of the basis braids for $l=1,2,\ldots,m$. 
 
%%%%%%%%%%%%%%%%%%%%%%%%%%%%%
\subsection{Triple points of $\mathcal{S}_m(a,b)$} \label{0203-1}
\begin{figure}
\begin{center}
  \includegraphics*{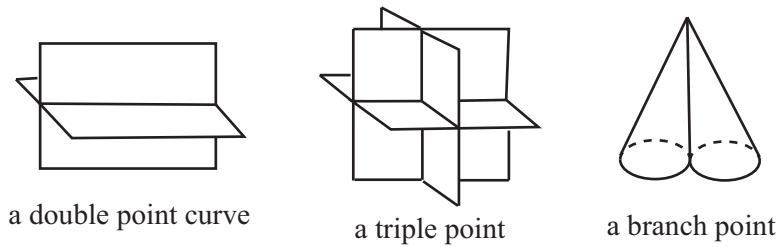}%%%%%%%%%%%%%%Fig.2.5%%%%%%%%%%%%%%%
\end{center}
  \caption{The singularity of a surface diagram.}
  \label{0215-1}
  \end{figure}

%%%%%%%%%%%%%%%%%%%%%%%%%%%%%%%%%5
Triple linking numbers are defined using triple points of a surface diagram. We will review a surface diagram of a surface link $S$ (see \cite{CKS}). For a projection $\pi \,:\, \mathbb{R}^4 \to \mathbb{R}^3$, the closure of the self-intersection set of $\pi(S)$ is called the singularity set. Let $\pi$ be a generic projection, i.e. the singularity set of the image $\pi(S)$ consists of double points, isolated triple points, and isolated branch points; see Fig. \ref{0215-1}. The closure of the singularity set forms a union of immersed arcs and loops, which we call double point curves. Triple points (resp. branch points) form the intersection points (resp. the end points) of the double point curves. A {\it surface diagram} of $S$ is the image $\pi(S)$ equipped with over/under information along each double point curve with respect to the projection direction. 

Throughout this paper, we consider the surface diagram of $\mathcal{S}_m(a,b)$ by the projection which projects $N(T)=I \times I \times T$ to $I \times T$ for an interval $I$, where we identify $N(T)$ with $I \times I \times T$ in such a way as follows. Since $T$ is the boundary of the standard solid torus in $\mathbb{R}^3 \times \{0\}$, the normal bundle of $T$ in $\mathbb{R}^3 \times \{0\}$ is a trivial bundle. We identify it with $I \times T$. Then we identify $N(T)$ with $I \times I \times T$, where the second $I$ is an interval in the fourth axis of $\mathbb{R}^4$. Perturbing $\mathcal{S}_m(a,b)$ if necessary, we can assume that this projection is generic. We call this surface diagram the surface diagram of $\mathcal{S}_m(a,b)$ in the form of a covering over the torus. 

Around a triple point $t$ in a surface diagram of an oriented surface link $S$, there are three types of sheets with respect to the projection direction. We call them the {\it top sheet}, {\it middle sheet} and {\it bottom sheet} from the higher one. We define that the {\it sign} of $t$ is $+1$ (resp. $-1$) if $t$ is a positive (resp. negative) triple point, and we call the triplet $(i, j, k)$ the {\it type} of $t$, where the top, middle, and bottom sheet are from the $i$th, $j$th, and $k$th component of $S$ respectively (\cite{CJKLS}); see Fig. \ref{0104-1}. The triple linking number $\mathrm{Tlk}_{i,j,k}(S)$ is the sum of the signs of all the triple points of type $(i,j,k)$ (\cite{CJKLS}). 

We can obtain the sign and type of each triple point of $\mathcal{S}_m(a,b)$ by studying the braid word transformation sequence from $ab$ to $ba$, as follows. Let us consider the surface obtained from $T$ by removing a sufficiently small open neighborhood of $\mathbf{m} \cup \mathbf{l}$, and denote it by $X$. The boundary of $X$ is a circle. Let us take two points $x$ and $x^\prime$ in $\partial X$ which divide $\partial X$ into two paths $\rho_0$ and $\rho_1$ from $x$ to $x^\prime$ such that $\rho_0$ (resp. $\rho_1$) is presented by $\mathbf{m} \cdot \mathbf{l}$ (resp. $\mathbf{l} \cdot \mathbf{m}$). Take an isotopy $\{\rho_u\}_{u \in [0,1]}$ of paths in $X$ which connects $\rho_0$ and $\rho_1$ such that $\rho_u \cap \rho_v =\{x, x^\prime \}$ for $u \neq v$. 

Put $S_u=\mathcal{S}_m(a,b) \cap p^{-1}(\rho_u)$ ($u \in [0,1]$). Since $\mathcal{S}_m(a,b)$ is a covering over $T$, $\{S_u\}_u$ is an isotopy of classical $m$-braids such that $S_0=ab$ and $S_1=ba$. Thus $\{S_u\}_u$ is presented by a braid word transformation sequence (see \cite{Kamada3}) from $ab$ to $ba$ related by a finite sequence of the following transformations: 
\begin{enumerate}[(1)]
\item
Insertion of the word $\sigma_i \sigma_i^{-1}$ or $\sigma_i^{-1} \sigma_i$, 
or deletion of the word $\sigma_i \sigma_i^{-1}$ or $\sigma_i^{-1} \sigma_i$, 
\item
Substitution of $\sigma_i \sigma_j$ for $\sigma_j \sigma_i$ where $|i-j|>1$,  
\item
Substitution of $\sigma_i \sigma_j \sigma_i$ for $\sigma_j \sigma_i \sigma_j$ where $|i-j|=1$.  
\end{enumerate}
Thus we have a braid word transformation sequence from $ab$ to $ba$. Conversely, when we have a braid word transformation sequence from $ab$ to $ba$, we can construct $\mathcal{S}_m(a,b)$. 

In the part of the surface diagram related by the transformation (1) or (2), the singularity set consists of disjoint double point curves; thus it does not contain triple points. For the transformation (3), we have the following lemma. The transformation (3) presents a Reidemeister move of type III. From now on we will call it a Reidemeister move of type III. 

\begin{lem} \label{0125-2}%%%%%%%%%%%%%%%%%%%%
Let us consider a braid word transformation sequence from $ab$ to $ba$ for commutative pure $m$-braids $a$ and $b$. We consider the surface diagram of $\mathcal{S}_m(a,b)$ in the form of a covering over $T$. Then the triple points are in the part of the surface diagram iff it is related by Reidemeister moves of type III (the transformation (3)), and the sign and type of each triple point are as follows. Let 
\begin{equation} \label{0125-1}
w \sigma_i \sigma_{j} \sigma_i w^\prime \rightarrow w \sigma_{j} \sigma_i \sigma_{j} w^\prime \  (|i-j|=1) 
\end{equation}
be a transformation in a braid word transformation sequence from $ab$ to $ba$, where $w$ and $w^\prime$ are some $m$-braids. The part of $\mathcal{S}_m(a,b)$ presented by this transformation has one triple point. The sign and type of the triple point is $+(\tau(j+1),  \tau(j), \tau(i))$ if $i<j$, and $-(\tau(i+1),  \tau(i), \tau(j))$ if $i>j$, where $\tau$ is the permutation associated with the $m$-braid $w^{-1}$. 
\end{lem}%%%%%%%%%%%%%%%%%%%%%

\Proof 
It suffices to determine the sign and type of the triple point presented by the transformation (\ref{0125-1}). Let us denote the triple point by $t$. 

\begin{figure}
\begin{center}
 \includegraphics*{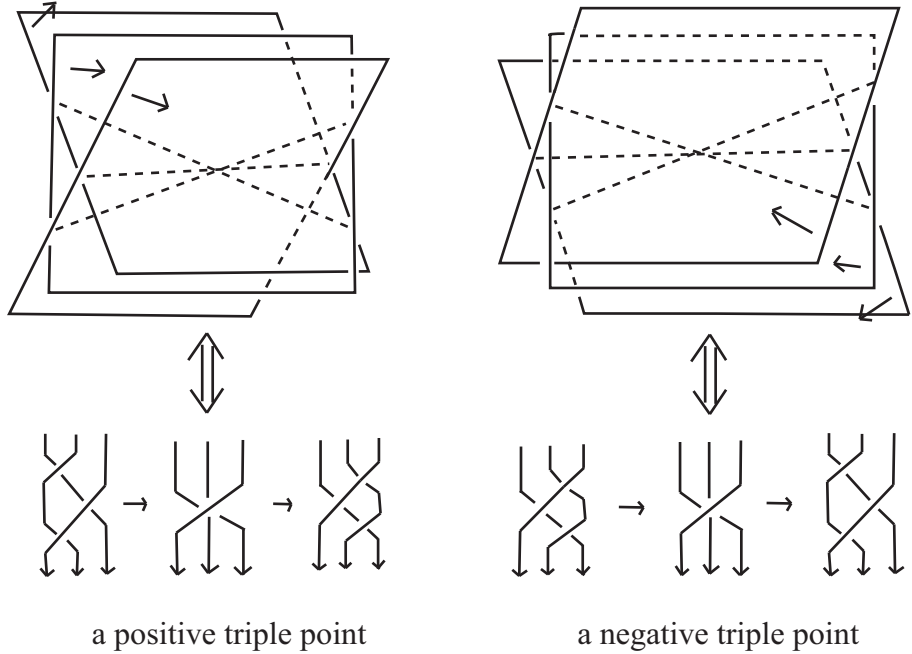}%%%%%%%%%%%%%%Fig.2.5%%%%%%%%%%%%%%%
\end{center}
 \caption{A positive triple point and a negative triple point. }
 \label{0122-1}
  \end{figure}

First we determine the sign of $t$. A positive triple point and a negative triple point as in Fig. \ref{0104-1} can be deformed as in Fig. \ref{0122-1}; thus a positive triple point is presented by a braid word transformation $\sigma_k \sigma_{k+1} \sigma_k \rightarrow \sigma_{k+1} \sigma_k \sigma_{k+1}$, and a negative triple point is presented by $\sigma_{k+1} \sigma_k \sigma_{k+1} \rightarrow \sigma_k \sigma_{k+1} \sigma_k$ (see also \cite[Proposition 4.43 (3)]{Carter-Saito}). Thus the sign of $t$ is positive (resp. negative) if $i<j$ (resp. $i>j$). 

Next we determine the type of $t$. The triple point $t$ is constructed by the $\tau(k)$th, $\tau(k+1)$th, and $\tau(k+2)$th strings of $m$-braids presented by (\ref{0125-1}). The top, middle, and bottom sheet are constructed by the $\tau(k+2)$th, $\tau(k+1)$th, and $\tau(k)$th strings. Since the $l$th string of $ab$ is from the $l$th component of $\mathcal{S}_m(a,b)$ for any $l$, so is the $l$th string of $m$-braids presented by (\ref{0125-1}). Thus the type of the triple point is $(\tau(k+2),  \tau(k+1), \tau(k))$. Thus the conclusion follows. 
\qed

%%%%%%%%%%%%%%%%%%%%%%%%%%%%%%5
\section{Triple linking numbers} \label{TripleLinking}
%%%%%%%%%%%%%%%%%%%%%%%%%%%%
In this section we prove Theorem \ref{linking}. 
\\

\noindent
%%%%%%%%%%%%%%%%%%%%%%%%%%%%%%%%%%%%%%%%%%55
{\it Proof of Theorem \ref{linking}.} 
%%%%%%%%%%%%%%%%%%%%%%%%%%%%%%%%%%%%%%%%%
Since triple linking numbers are determined from three components of a surface link, it is sufficient to consider three components of $\mathcal{S}_m(a,b)$. Since $\mathcal{S}_m(a,b)$ is in the form of a covering over $T$, any three components of $\mathcal{S}_m(a,b)$ form a torus-covering $T^2$-link. Thus it suffices to show for the case when $m=3$. From now on we consider two commutative pure $3$-braids $a$ and $b$. 

We determine the presentations of $a$ and $b$ as follows. It is known (see \cite{Hansen}) that the pure $3$-braid group $P_3$ is presented by generators $A_{12}$, $A_{23}$, $A_{13}$, and the relations
\begin{eqnarray*}
&& A_{12}^{-1} A_{23} A_{12}=A_{13} A_{23} A_{13}^{-1}, \\
&& A_{12}^{-1} A_{13} A_{12}=A_{13} A_{23} A_{13} A_{23}^{-1} A_{13}^{-1}, 
\end{eqnarray*}
where 
\[
A_{12}=\sigma_1^2, \ A_{23}=\sigma_2^2, \ A_{13}=\sigma_2 \sigma_1^2 \sigma_2^{-1}. 
\]
Put $A=A_{12}$, $B=A_{23}$ and $C=A_{13}$. Then 
\[
P_3=\langle A,B,C \mid A^{-1}BA=CBC^{-1}, A^{-1} CA=CBCB^{-1}C^{-1} \rangle. 
\]
Put $\Delta=ACB$; note that $\Delta$ is a full twist. By eliminating $C$ by $C=A^{-1} \Delta B^{-1}$, we have 
\begin{eqnarray*}
P_3 &=& \langle A,B,\Delta \mid A \Delta=\Delta A, B \Delta=\Delta B \rangle \\
     &=& Z \oplus U, 
\end{eqnarray*}
where $Z$ is an infinite cyclic group generated by $\Delta$, and $U$ is the free group with two generators $A$ and $B$; note that since the center of $U$ is trivial, the center of $P_3$ is $Z$. Thus we can see that any element of $P_3$ is written uniquely by $\Delta^l u$ for an integer $l$ and $u \in U$. Let $a=\Delta^{k_1} u$, $b=\Delta^{k_2} v$ (for integers $k_1$, $k_2$ and $u, v \in U$) be the presentations. Since $Z$ is the center of $P_3$, we can see that $a$ and $b$ commute iff $u$ and $v$ commute. Since $U$ is a free group with two generators, $u$ and $v$ commute iff $u=w^{l_1}$ and $v=w^{l_2}$ for some $w \in U$ and integers $l_1$ and $l_2$. Since $a$ and $b$ are commutative pure $3$-braids, $a=\Delta^{k_1} w^{l_1}$ and $b=\Delta^{k_2} w^{l_2}$. 

We have a braid word transformation sequence as follows: 
\begin{eqnarray*}%%%%%%%%%%%%%%%%%%%%%
ab &=& \Delta^{k_1} (w^{l_1} \Delta^{k_2}) w^{l_2} \\
&\rightarrow& \cdots \rightarrow  
\Delta^{k_1} (\Delta^{k_2} w^{l_1}) w^{l_2} =\Delta^{k_2} (\Delta^{k_1} w^{l_2}) w^{l_1} \\
&\rightarrow & \cdots \rightarrow 
\Delta^{k_2} (w^{l_2} \Delta^{k_1}) w^{l_1} =ba, 
\end{eqnarray*}%%%%%%%%%%%%%%%%%%%%%%%%%5
where a transformation from $w \Delta$ to $\Delta w$ (resp. from $\Delta w$ to $w \Delta$) is applied $l_1 k_2$ (resp. $k_1 l_2$) times. Since $\Delta$ and $w$ are pure braids, the triple linking numbers in the surface concerning a transformation from $w \Delta$ to $\Delta w$ (resp. from $\Delta w$ to $w \Delta$) have the same value with those of $\mathcal{S}_m(w, \Delta)$ (resp. $\mathcal{S}_m(\Delta, w)$). Hence 
\begin{eqnarray} \label{0202-1}
&& \mathrm{Tlk}_{i,j,k}(\mathcal{S}_3(a,b))\\
&=& l_1 k_2 \mathrm{Tlk}_{i,j,k}(\mathcal{S}_3(w,\Delta))
+ k_1 l_2 \mathrm{Tlk}_{i,j,k}(\mathcal{S}_3(\Delta, w)). \nonumber
\end{eqnarray}

It is known \cite{CJKLS} that for any oriented surface link $S$, $\mathrm{Tlk}_{k,j,i}(S)=-\mathrm{Tlk}_{i,j,k}(S)$ if $i$, $j$, $k$ are mutually distinct, and $\mathrm{Tlk}_{i,j,k}(S)=0$ otherwise. Hence it suffices to calculate $\mathrm{Tlk}_{i,j,k}(\mathcal{S}_3(a,b))$ for $(i,j,k)=(1,2,3)$, $(2,3,1)$ and $(3,1,2)$. By Lemma \ref{0126-2}, the triple linking numbers are presented by the exterior product as follows: 
\begin{equation*} \label{0126-3}%%%%%%%%%%%%%%%%%%%
\left(
\begin{array}{c}
  \mathrm{Tlk}_{1,2,3}(\mathcal{S}_3(c, \Delta)) \\
  \mathrm{Tlk}_{2,3,1}(\mathcal{S}_3(c, \Delta)) \\
  \mathrm{Tlk}_{3,1,2}(\mathcal{S}_3(c, \Delta)) 
\end{array}
\right)
=-\left(
\begin{array}{c}
  \mathrm{Lk}_{3,1}(\hat{c}) \\
  \mathrm{Lk}_{1,2}(\hat{c}) \\
  \mathrm{Lk}_{2,3}(\hat{c}) 
\end{array}
\right)
\times
\left(
\begin{array}{c}
 1 \\
 1 \\
 1
\end{array}
\right)
\end{equation*}
where $c=A$ or $B$; note that $(\mathrm{Lk}_{3,1}(\hat{\Delta}), \mathrm{Lk}_{1,2}(\hat{\Delta}), \mathrm{Lk}_{2,3}(\hat{\Delta}))=(1,1,1)$. Since any element of $U$ is presented by $A$ and $B$, it follows that this equation holds for any $c \in U$. Put $\mathbf{w}=^t (
  \mathrm{Lk}_{3,1}(\hat{w}) \ 
  \mathrm{Lk}_{1,2}(\hat{w}) \ 
  \mathrm{Lk}_{2,3}(\hat{w}) 
)$. 
Since $w \in U$, by (\ref{0202-1}) and Lemma \ref{0126-2}, 
\begin{eqnarray*}
\left(
\begin{array}{c}
  \mathrm{Tlk}_{1,2,3}(\mathcal{S}_3(a,b)) \\
  \mathrm{Tlk}_{2,3,1}(\mathcal{S}_3(a,b)) \\
  \mathrm{Tlk}_{3,1,2}(\mathcal{S}_3(a,b)) 
\end{array}
\right)
&=& -l_1k_2 \mathbf{w} \times 
\left(
\begin{array}{c}
 1 \\
 1 \\
 1
\end{array}
\right)
- k_1 l_2\left(
\begin{array}{c}
 1 \\
 1 \\
 1
\end{array}
\right)
\times \mathbf{w}
\\
&=& -( 
l_1 \mathbf{w} \times 
k_2 \left(
\begin{array}{c}
 1 \\
 1 \\
 1
\end{array}
\right)
+
k_1 \left(
\begin{array}{c}
 1 \\
 1 \\
 1
\end{array}
\right) 
\times l_2 \mathbf{w} )
\\
&=&
-( k_1 \left(
\begin{array}{c}
 1 \\
 1 \\
 1
\end{array}
\right) +l_1 \mathbf{w}
)
\times
( k_2 \left(
\begin{array}{c}
 1 \\
 1 \\
 1
\end{array}
\right) +l_2 \mathbf{w}
)
\\%
&=& -\left(
\begin{array}{c}
  \mathrm{Lk}_{3,1}(\hat{a}) \\
  \mathrm{Lk}_{1,2}(\hat{a}) \\
  \mathrm{Lk}_{2,3}(\hat{a})
\end{array}
\right)
\times
\left(
\begin{array}{c}
   \mathrm{Lk}_{3,1}(\hat{b}) \\
  \mathrm{Lk}_{1,2}(\hat{b}) \\
  \mathrm{Lk}_{2,3}(\hat{b})
\end{array}
\right). 
\end{eqnarray*}
Thus 
$\mathrm{Tlk}_{i,j,k}(\mathcal{S}_m(a,b))=-\mathrm{Lk}_{i,j}(\hat{a}) \mathrm{Lk}_{j,k}(\hat{b})+\mathrm{Lk}_{i,j}(\hat{b}) \mathrm{Lk}_{j,k}(\hat{a})$. 
\qed
\\

We can check that for the $T^2$-link $S$ of Lemma \ref{0126-2}, $\mathrm{Tlk}_{k,j,i}(S)=-\mathrm{Tlk}_{i,j,k}(S)$ if $i$, $j$, $k$ are mutually distinct, and otherwise $\mathrm{Tlk}_{i,j,k}(S)=0$. Thus we can see that this property holds for $\mathcal{S}_m(a,b)$ without referring to \cite{CJKLS}. 

\begin{lem} \label{0126-2}%%%%%%%%%%%%%%%%%
Put $A=\sigma_1^2$ and $B=\sigma_2^2$ as in the proof of Theorem \ref{linking}. Then 
\begin{eqnarray}
&& \mathrm{Tlk}_{i,j,k}(\mathcal{S}_3(c,\Delta))=-\mathrm{Lk}_{i,j}(\hat{c}) + \mathrm{Lk}_{j,k}(\hat{c}) \label{0124-1}\\
&& \mathrm{Tlk}_{i,j,k}(\mathcal{S}_3(\Delta,c))=-\mathrm{Tlk}_{i,j,k}(\mathcal{S}_3(c,\Delta)) \label{0124-3} 
\end{eqnarray}
where $c=A$ or $B$. 
\end{lem}

\Proof
We show (\ref{0124-1}) for the case when $c=A=\sigma_1^2$. Since $\Delta$ is a full twist, we assume that $\Delta$ has a presentation $\Delta=\sigma_1 \sigma_2 \sigma_1 \sigma_2 \sigma_1 \sigma_2$. We have the following braid word transformation sequence: 
\begin{eqnarray}
\sigma_1 \Delta &=& \sigma_1 (\sigma_1 \sigma_2 \sigma_1) \sigma_2 \sigma_1 \sigma_2 \label{0124-5}\\
& \rightarrow & \sigma_1 (\sigma_2 \sigma_1 \sigma_2) \sigma_2 \sigma_1 \sigma_2
=\sigma_1 \sigma_2 \sigma_1 \sigma_2 (\sigma_2 \sigma_1 \sigma_2) \nonumber \\
& \rightarrow & \sigma_1 \sigma_2 \sigma_1 \sigma_2 (\sigma_1 \sigma_2 \sigma_1) 
=\Delta \sigma_1. \nonumber
\end{eqnarray}
Since a Reidemeister move of type III is applied twice, there are two triple points in the surface diagram presented by (\ref{0124-5}): the first one is positive, and the second one is negative by Lemma \ref{0125-2}. Since $A=\sigma_1^2$, repeating (\ref{0124-5}) twice, we have the braid word transformation sequence 
\begin{equation} \label{0126-1}%%%%%%%%%%%%%%%%%%%
A \Delta=\sigma_1^2 \Delta \rightarrow \sigma_1 D \rightarrow \sigma_1 \Delta \sigma_1 
\rightarrow D \sigma_1 \rightarrow \Delta \sigma_1^2=\Delta A, 
\end{equation}%%%%%%%%%%%%%%%%%%%%%%%%%%%%%%%%
where $D=\sigma_1 \sigma_2 \sigma_1 \sigma_2^2 \sigma_1 \sigma_2$. Thus there are four triple points. We denote them by $t_1, \ldots, t_4$. The sign of $t_1$ and $t_3$ (resp. $t_2$ and $t_4$) is positive (resp. negative) by Lemma \ref{0125-2}. 

Next we determine the type of $t_l$ ($l=1,\ldots,4$). For the $m$-braid group $B_m$ and the symmetric group $S_m$ of degree $m$, let $\phi \,:\, B_m \rightarrow S_m$ be a homomorphism which maps an $m$-braid $c$ to a permutation associated with $c$, i.e. $\phi(\sigma_l)=(l \ l+1)$ for $l=1,2,\ldots, m-1$. Let $w_l$ ($l=1, \ldots,4$) be the braid $w$ of the transformation (\ref{0125-1}) presenting $t_l$ (see Lemma \ref{0125-2}). Then $w_1=\sigma_1^2$, $w_2=\sigma_1^2 \sigma_2 \sigma_1 \sigma_2$, $w_3=\sigma_1^{-1} w_1$, and $w_4=\sigma_1^{-1} w_2$. Put $\tau_i=\phi(w_i^{-1}) \in S_3$ ($i=1,\ldots, 4$). Since $\tau_1=e$ (resp. $\tau_2=(1\ 3)$), the type of $t_1$ (resp. $t_2$) is $(\tau_1(3), \tau_1(2), \tau_1(1))=(3,2,1)$ (resp. $(\tau_2(3), \tau_2(2), \tau_2(1))=(1,2,3)$). Since $w_3^{-1}=w_1^{-1} \sigma_1$ (resp. $w_4^{-1}=w_2^{-1} \sigma_1$), we can obtain $\tau_3$ (resp. $\tau_4$) by applying the action of $\phi(\sigma_1)=(1\ 2)$ on $\tau_1$ (resp. $\tau_2$); thus the type of $t_3$ (resp. $t_4$) is $(3,1,2)$ (resp. $(2,1,3)$). 

Thus the sign and type of all the triple points are 
\[
+(3,2,1), -(1,2,3), +(3,1,2), -(2,1,3), 
\]
and hence 
\[
\mathrm{Tlk}_{i,j,k}(\mathcal{S}_3(A, \Delta))=\begin{cases} +1 & \mathrm{if} \ (i,j,k)=(3,2,1) \ \mathrm{or} \ (3,1,2), \\
                                            -1& \mathrm{if} \ (i,j,k)=(1,2,3) \ \mathrm{or} \ (2,1,3), \\
                                              0 & \mathrm{otherwise}. 
\end{cases}
\]
Since $A=\sigma_1^2 \in P_3$, $\mathrm{Lk}_{i,j}(\hat{A})=+1$ if $\{i,j\}=\{1,2\}$ and otherwise zero. Hence $\mathrm{Tlk}_{i,j,k}(\mathcal{S}_3(A, \Delta))=-\mathrm{Lk}_{i,j}(\hat{A}) + \mathrm{Lk}_{j,k}(\hat{A})$ for any $i,j,k$.  

Next we show (\ref{0124-3}) for the case when $c=A$. We have the braid word transformation sequence from $\Delta A$ to $A \Delta$ which is the reversed sequence of (\ref{0126-1}); thus the sign of each triple point is the reverse of the sign of $t_l$, and the type is the same with that of $t_l$ ($l=1,\ldots ,4$). Hence $\mathrm{Tlk}_{i,j,k}(\mathcal{S}_3(A, \Delta))=-\mathrm{Tlk}_{i,j,k}(\mathcal{S}_3(A, \Delta))$ for any $i,j,k$. 
 
Since $B=\sigma_2^2$ (resp. $\Delta$) is obtained from $A=\sigma_1^2$ (resp. $\Delta$) by regarding the $l$th string as the $(3-l)$th string for $l=1,2,3$, we can see that (\ref{0124-1}) and (\ref{0124-3}) also hold for the case when $c=B$. 
\qed
   
%%%%%%%%%%%%%%%%%%%%%%%%%%%%%%%%%%%%
\section{Triple point numbers} \label{0104-3}
%%%%%%%%%%%%%%%%%%%%%%%%%%%%%%%%%
  
In this section, we show Theorem \ref{thm-triple-pt} by using Theorem \ref{linking}. Further, we show Theorem \ref{example} by using Theorem \ref{thm-triple-pt}. 

 It is known \cite{CKSS01} (see also \cite{CKSS02}) that any oriented surface link is bordant to the split union of oriented \lq\lq necklaces", and \cite{CKSS02} any surface link is unorientedly bordant to the split union of necklaces and connected sums of standard projective planes; see also \cite{Sanderson}. The triple point number of a necklace is easily determined to be $4n$, as follows. A \lq\lq necklace" is a surface link introduced in \cite{CKSS02} with non-trivial triple linking (see also \cite{CJKLS, CKS, CKSS01}). An oriented necklace is a \lq\lq Hopf 2-link without or with beads" (\cite{CKSS01}). A necklace consists of two components homeomorphic to two tori or two Klein bottles (called a \lq\lq strand") and $n$ components homeomorphic to 2-spheres (called \lq\lq beads"); see \cite{CKSS02}. By estimating the lower bound of the triple point number by using the triple linking numbers (see \cite{CKSS01, CKSS02}, see \cite{Satoh2} for the triple linking numbers of a non-orientable surface link), we can see that the triple point number concerning each bead is $4$, and as a whole the triple point number is $4n$. 
 
%%%%%%%%%%%%%%%%%%%%%%%%%%%%%%%%%%%%%%%%%%%%5
\subsection{Proof of Theorem \ref{thm-triple-pt}}
%%%%%%%%%%%%%%%%%%%%%%%%%%%%%%%%%%%%%%%%%%
 Put $S=\mathcal{S}_m(b, \Delta^n)$. By definition, $t(S) \geq \sum_{i \neq j, j \neq k}|\mathrm{Tlk}_{i,j,k}(S)|$. We show that $\sum_{i \neq j, j \neq k}|\mathrm{Tlk}_{i,j,k}(S)|=4n(\mu (m-2)-\nu)$, where $\mu=\sum_{i<j}|\mathrm{Lk}_{i,j}(\hat{b})|$, and $\nu=\sum_{i<j<k} (\nu_{i,j,k}+\nu_{j,k,i}+\nu_{k,i,j})$, where $\nu_{i,j,k}=\min_{i,j,k} \{ |\mathrm{Lk}_{i,j}(\hat{b})|, |\mathrm{Lk}_{j,k}(\hat{b})| \}$ if $\mathrm{Lk}_{i,j}(\hat{b}) \mathrm{Lk}_{j,k}(\hat{b})>0$ and otherwise zero. Since $\mathrm{Lk}_{i,j}(\hat{\Delta^n})=n$ for any $i \neq j$, 
\[
\mathrm{Tlk}_{i,j,k}(S)=-n(\mathrm{Lk}_{i,j}(\hat{b})-\mathrm{Lk}_{j,k}(\hat{b}))
\]
for any $i,j,k$ with $i \neq j$ and $j \neq k$ by Theorem \ref{linking}. 

It suffices to show that 
\[
\sum_{i \neq j, j \neq k} |\mathrm{Lk}_{i,j}(\hat{b})-\mathrm{Lk}_{j,k}(\hat{b})|=4(\mu (m-2)-\nu). 
\]
By definition of $\nu_{i,j,k}$, 
\[
|\mathrm{Lk}_{i,j}(\hat{b})-\mathrm{Lk}_{j,k}(\hat{b})|=|\mathrm{Lk}_{i,j}(\hat{b})|+|\mathrm{Lk}_{j,k}(\hat{b})|-2\nu_{i,j,k}. 
\]
We calculate the sum of $|\mathrm{Lk}_{i,j}(\hat{b})|+|\mathrm{Lk}_{j,k}(\hat{b})|$ as follows. 
\begin{eqnarray*}
\sum_{i\neq j, j \neq k} |\mathrm{Lk}_{i,j}(\hat{b})|+|\mathrm{Lk}_{j,k}(\hat{b})|
&=& 2 \sum_{i\neq j, j \neq k}|\mathrm{Lk}_{i,j}(\hat{b})| \\
&=&  2(m-1)\sum_{i\neq j} |\mathrm{Lk}_{i,j}(\hat{b})| \\
&=& 4 \mu (m-1). 
\end{eqnarray*}
We calculate the sum of $\nu_{i,j,k}$ as follows. 
\begin{eqnarray*}
\sum_{i\neq j, j \neq k} \nu_{i,j,k} 
&=& \sum_{i\neq j, j \neq k, k \neq i} \nu_{i,j,k} +\sum_{i\neq j} |\mathrm{Lk}_{i,j}(\hat{b})| \\
&=& 2 \sum_{i<j<k}(\nu_{i,j,k} + \nu_{j,k,i} + \nu_{k,i,j}) +2 \mu \\
&=& 2(\nu+\mu). 
\end{eqnarray*}
Thus 
\begin{eqnarray*}
\sum_{i \neq j, j \neq k} |\mathrm{Lk}_{i,j}(\hat{b})-\mathrm{Lk}_{j,k}(\hat{b})|
&=& 4 \mu (m-1)-2 \cdot 2(\nu+\mu) \\
&=& 4(\mu(m-2)-\nu), 
\end{eqnarray*}
and hence $t(S) \geq 4n(\mu(m-2)-\nu)$. 
\qed

 %%%%%%%%%%%%%%%%%%%%%%%%%%%%%%%%%%%%%%%%%%%%5
\subsection{Proof of Theorem \ref{example}}
%%%%%%%%%%%%%%%%%%%%%%%%%%%%%%%%%%%%%%%%%%
Put $S=\mathcal{S}_m(b, \Delta^n)$. By Theorem \ref{thm-triple-pt}, $t(S) \geq 4n(\mu(m-2)-\nu)$. Since $b$ consists of $\sigma_i^{2(-1)^{i}}$ ($i=1,2,\ldots, m-1$), $\mathrm{Lk}_{i,j}(\hat{b})$ is positive (resp. negative) if $j=i+1$ and $i$ is even (resp. $j=i+1$ and $i$ is odd), and otherwise zero. Thus $\nu_{i,j,k}=0$ for any $i,j,k$, and it follows that $\nu=0$. Hence $t(S) \geq 4n \mu (m-2)$, where $\mu=\sum_{i<j}|\mathrm{Lk}_{i,j}(\hat{b})|$. 

Now we show that the surface diagram of $S$ in the form of a covering over $T$ (see Section \ref{0203-1}) has $4n \mu (m-2)$ triple points. By Lemma \ref{0125-2}, it suffices to show that there is a braid word transformation sequence from $b \Delta^n$ to $\Delta^n b$ related $4n \mu (m-2)$ times by a Reidemeister move of type III. 

By the presentation of $b$, $b$ consists of $2 \mu$ standard generators; thus it suffices to show that there is a braid word transformation sequence from $\sigma_i \Delta$ to $\Delta \sigma_i$ related $2(m-2)$ times by a Reidemeister move of type III, as follows. We assume that a full twist $\Delta$ has the following presentation: 
\[
\Delta=(\sigma_1 \sigma_2 \cdots \sigma_{m-1})^m.
\]
The transformation sequence from $\sigma_i \Delta$ to $\Delta \sigma_i$ consists of the following parts: 
\begin{enumerate}[(i)]
\item 
A transformation from $\sigma_j (\sigma_1 \sigma_2 \cdots \sigma_{m-1})$ to $(\sigma_1 \sigma_2 \cdots \sigma_{m-1}) \sigma_{j-1}$, where $2 \leq j \leq m-1$.  

\item 
A transformation from $\sigma_1 (\sigma_1 \sigma_2 \cdots \sigma_{m-1})^2$ to $(\sigma_1 \sigma_2 \cdots \sigma_{m-1})^2 \sigma_{m-1}$.  
\end{enumerate}
We can take a transformation sequence (i) (resp. (ii)) related once (resp. $(m-2)$ times) by a Reidemeister move of type III by Lemma \ref{0201-1}. The transformation sequence from $\sigma_i \Delta$ to $\Delta \sigma_i$ is related $(i-1)$ times by the transformation (i), and then once by the transformation (ii), and then $(m-i-1)$ times by the transformation (i). Thus the total applied number is $(i-1)+(m-2)+(m-i-1)=2(m-2)$. Thus $t(S)=4n(m-2)(\sum_{i<j}|\mathrm{Lk}_{i,j}(\hat{b})|)$, and this is realized by the surface diagram in the form of a covering over the torus. 
\qed

\begin{lem} \label{0201-1}%%%%%%%%%%%%%%%%%%%
We can take a transformation sequence (i) (resp. (ii)) as in the proof of Theorem \ref{example} which is related once (resp. $(m-2)$ times) by a Reidemeister move of type III. 
\end{lem}%%%%%%%%%%%%%%%%%%%%%%%%%%%%%%%
\Proof 
Let us regard two braids as equivalent if they are related by transformations other than Reidemeister moves of type III. Then we can take a transformation sequence (i) as follows. 
\begin{eqnarray*}
 \sigma_j (\sigma_1 \sigma_2 \cdots \sigma_{m-1}) 
&\sim& 
\sigma_1 \sigma_2 \cdots \sigma_{j-2} (\sigma_{j} \sigma_{j-1} \sigma_j) \sigma_{j+1} \cdots \sigma_{m-1} \\ 
&\rightarrow& \sigma_1 \sigma_2 \cdots \sigma_{j-2} 
(\sigma_{j-1} \sigma_{j} \sigma_{j-1}) \sigma_{j+1} \cdots \sigma_{m-1} \\
&\sim& (\sigma_1 \sigma_2 \cdots \sigma_{m-1}) \sigma_{j-1}. 
\end{eqnarray*}
Thus the transformation sequence (i) is related once by a Reidemeister move of type III. 
 
We can take a transformation sequence (ii) as follows. 
\begin{eqnarray*}
&& \sigma_1 (\sigma_1 \sigma_2 \cdots \sigma_{m-1})^2 \\
&\sim& 
\sigma_1 (\sigma_1) (\sigma_2 \sigma_1) (\sigma_3 \sigma_2) \cdots (\sigma_{j} \sigma_{j-1}) \cdots 
(\sigma_{m-1} \sigma_{m-2}) \sigma_{m-1} \\
&=&
\sigma_1 (\sigma_1 \sigma_2 \sigma_1) \sigma_3 \sigma_2 \cdots 
\sigma_{m-1} \sigma_{m-2} \sigma_{m-1} \\ 
&\rightarrow& 
\sigma_1 (\sigma_2 \sigma_1 \sigma_2) \sigma_3 \sigma_2 \cdots  
\sigma_{m-1} \sigma_{m-2} \sigma_{m-1} \\
&= &
\sigma_1 \sigma_2 \sigma_1 (\sigma_2 \sigma_3 \sigma_2) \cdots  
\sigma_{m-1} \sigma_{m-2} \sigma_{m-1} \\ 
&\rightarrow& 
\sigma_1 \sigma_2 \sigma_1 (\sigma_3 \sigma_2 \sigma_3) \cdots 
\sigma_{m-1} \sigma_{m-2} \sigma_{m-1} 
\\ 
&\rightarrow& \cdots \\
& \rightarrow& 
\sigma_1 \sigma_2 \sigma_1 \sigma_3 \sigma_2  \cdots  
(\sigma_{m-2} \sigma_{m-1} \sigma_{m-2}) \sigma_{m-1} 
\\ 
&\rightarrow& 
\sigma_1 \sigma_2 \sigma_1 \sigma_3 \sigma_2  \cdots 
 (\sigma_{m-1} \sigma_{m-2} \sigma_{m-1}) \sigma_{m-1} \\ 
&\sim& (\sigma_1 \sigma_2 \cdots \sigma_{m-1})^2 \sigma_{m-1}. 
\end{eqnarray*}
Thus the transformation sequence (ii) is related $(m-2)$ times by Reidemeister moves of type III. 
\qed
\\

    \section*{Acknowledgements} 
The author would like to thank Professors Shin Satoh and Seiichi Kamada for their helpful comments. The author is supported by GCOE, Kyoto University.


\begin{thebibliography}{0}

\bibitem{CJKLS}
J. S. Carter, D. Jelsovsky, S. Kamada, L. Langford and M. Saito, 
{\it Quandle cohomology and state-sum invariants of knotted curves and surfaces}, 
Trans. Amer. Math. Soc. {\bf 355} (2003), 3947--3989.

\bibitem{CKS01}
J. S. Carter, S. Kamada, and M. Saito, 
{\it Geometric interpretations of quandle homology and cocycle knot invariant}, 
J. Knot Theory Ramificartions {\bf 10} (2001), 345--358. 

\bibitem{CKS}
J. S. Carter, S. Kamada, M. Saito, 
{\it Surfaces in 4-space}, Encyclopaedia of Mathematical Sciences 142, 
Low-Dimensional Topology III (Springer-Verlag, Berlin, 2004). 

\bibitem{CKSS01}
J. S. Carter, S. Kamada, M. Saito and S. Satoh, 
{\it A theorem of Sanderson on link bordisms in dimension $4$}, 
Algebraic and Geometric Topology {\bf 1} (2001), 299--310. 

\bibitem{CKSS02}
J. S. Carter, S. Kamada, M. Saito and S. Satoh, 
{\it Bordism of unoriented surfaces in $4$-space}, 
Michigan Math. {\bf 50} (2002), 575--591. 

\bibitem{COS}
J. S. Carter, K. Oshiro and M. Saito, 
{\it Symmetric extensions of dihedral quandles and triple points of non-orientable surfaces}, 
Topology Appl. {\bf 157} (2010), 857--869.


\bibitem{Carter-Saito}
 J. S. Carter and M. Saito, {\it Knotted surfaces and their diagrams}, Mathematical Surveys 
 and Monographs 55 (Amer. Math. Soc., 1998).
  
 
\bibitem{Hansen}
V. L. Hansen, {\it Braids and Coverings: Selected Topics}, London Mathematical Society 
Student Texts, Vol. 18 (Cambridge University Press, 1989). 

\bibitem{Hatakenaka}
E. Hatakenaka, {\it An estimate of the triple point numbers
of surface-knots by quandle cocycle invariants}, 
 Topology Appl. {\bf 139}  (2004), 129--144.

\bibitem{Iwakiri}
M. Iwakiri, {\it Triple point cancelling numbers of surface links and quandle cocycle invariants}, 
Topology Appl. {\bf 153} (2006) 2815--2822. 

 
\bibitem{Kamada1}
S. Kamada, {\it A characterization of groups of closed orientable surfaces in 4-space}, 
Topology {\bf 33} (1994), 113-122.

 
 \bibitem {Kamada3}
 S. Kamada, {\it Braid and Knot Theory in Dimension Four},
Math. Surveys and Monographs 95 (Amer. Math. Soc., 2002). 
  

\bibitem{Kamada-Oshiro}
S. Kamada and K. Oshiro, 
{\it Homology groups of symmetric quandles and cocycle invariants of links and 
surface links}, Trans. Amer. Math. Soc. {\bf 362} (2010), 5501--5527.

 \bibitem{N}
I. Nakamura, {\it Surface links which are coverings over the standard torus}, 
Algebraic and Geometric Topology {\bf 11} (2011), 1497--1540. 

\bibitem{Oshiro10}
K. Oshiro, {\it Triple point numbers of surface-links and symmetric
quandle cocycle invariants}, 
Algebraic and Geometric Topology {\bf 10} (2010), 853--865. 
 

\bibitem{Rolfsen}
D. Rolfsen, {\it Knots and Links}, (Publish or Perish Press, Berkley, 1976). 


\bibitem{Sanderson}
B. J. Sanderson, {\it Bordism of links in codimension 2}, 
J. London Math. Soc. {\bf 35} (1987), no. 2, 
367--376. 


\bibitem{Satoh}
S. Satoh, 
{\it Minimal triple numbers of some non-orientable surface-links}, 
Pacific J. Math. {\bf 197} (2001), 213--221. 

\bibitem{Satoh2}
S. Satoh, {\it Triple point invariants of non-orientable surface-links}, 
Topology Appl. {\bf 121} (2002), 207--218. 

\bibitem{SatShi04}
S. Satoh and A. Shima, 
{\it The $2$-twist-spun trefoil has the triple point number four}, 
Trans. Amer. Math. Soc. {\bf 356} (2004), 1007--1024. 

\bibitem{SatShi05}
S. Satoh and A. Shima, 
{\it Triple point numbers and quandle cocycle invariants of knotted surfaces 
in 4-space}, New Zealand J. Math. {\bf 34} (2005), no. 1, 71--79. 

\end{thebibliography}
\end{document}